\pdfoutput=1
\documentclass[%
 aip,
 jmp,%
 amsmath,amssymb,
 reprint,%
]{revtex4-1}

\usepackage{graphicx}
\usepackage{dcolumn}
\usepackage{bm}
\usepackage{blindtext}
\usepackage{mathtools}
\usepackage{amsmath}
\usepackage[ruled]{algorithm2e}
\usepackage{algorithmic}
\usepackage{paracol}
\usepackage{capt-of}
\def\x{{\mathbf x}}
\def\p{{\mathbf p}}

\renewcommand\appendixname{Derivation of the proposed dynamical system model:}

\begin{document}


\title[Dynamical Global Optimization]{
Global Optimization based on Growth Transform Dynamical Model }

\author{Oindrila Chatterjee}
\author{Shantanu Chakrabartty}
\email{shantanu@wustl.edu}
\affiliation{ 
Department of Electrical and Systems Engineering, Washington University in St. Louis,\newline Missouri 63130, USA
}%

\date{\today}

\begin{abstract} 
Conservation principles like conservation of charge or energy provide a natural way to couple and constrain different physical variables. In this letter, we propose a dynamical system model that exploits these constraints for solving non-convex global optimization problems. Unlike the traditional simulated annealing or quantum annealing based global optimization techniques, the proposed method optimizes a target objective function by continuously evolving a driver functional over a conservation manifold using a generalized variant of growth transformations. As a result, the driver functional converges to a Dirac-delta function that is centered at the global optimum of the target objective function. We provide an outline of the proof of convergence for the dynamical system model and we demonstrate the application of the model for implementing linear-time and constant-time decentralized sorting algorithms.
\end{abstract}

\keywords{Global optimization, Analog computing, Growth transforms, Dynamical systems, Continuous-time systems}
\maketitle

\begin{figure}
\begin{center}
\includegraphics[page=1,scale=0.35,trim=4 4 4 4,clip]{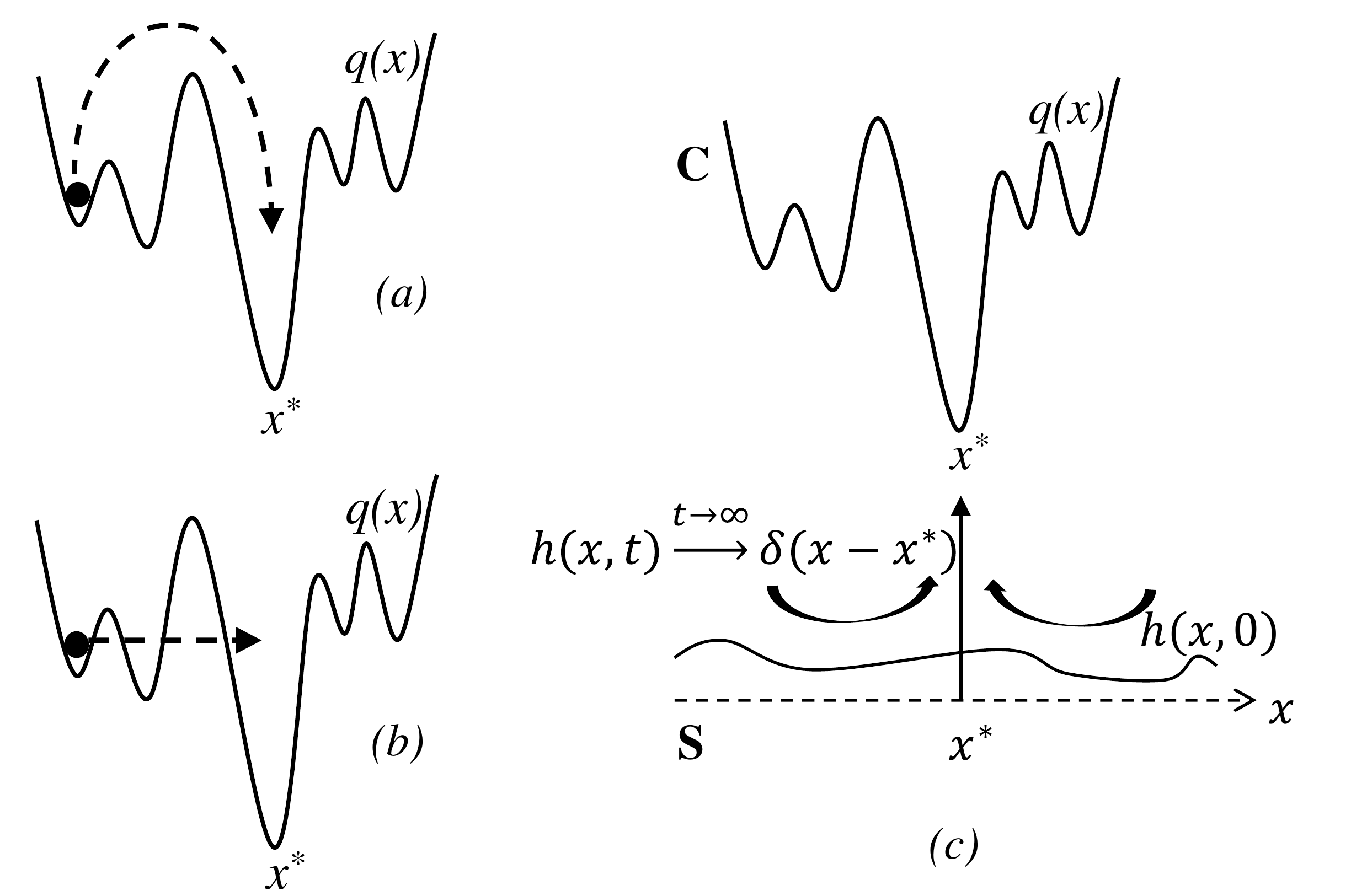}
\end{center}
\vspace{-0.5cm}
\caption{Illustration of different annealing principles for optimizing a target objective function $q(\x)$ with multiple local minima but only one global minimum $\x^{*}$: (a)the process of surmounting the energy-barriers in simulated annealing; (b) the process of quantum tunneling through the barriers in quantum annealing; and (c) the proposed approach where a driver function $h(\x,t)$ evolves under the influence of $q(\x)$ to a Dirac-delta function centered at the global minimum $\x^{*}$ . }
\label{fig_globalopt}
\end{figure}

Naturally occurring systems generally obey two physical principles: (a) an isolated system converges 
to a set of equilibrium states (referred to as eigenstates) that correspond to the lowest energy\cite{feynman2013feynman}, and (b) the dynamics 
of the system evolves in a manner that some physical quantities (for example energy, charge, mass or 
momentum) are conserved~\cite{feynman1967character}. In literature these two principles have independently served as the basis for designing
analog computing systems~\cite{ vergis1986complexity}, analog neural networks \cite{hopfield1985neural} and Lyapunov networks \cite{hopfield1995rapid}. 
These two physical principles have also formed the basis for annealing algorithms that have been used for optimizing complex, non-convex functions with multiple
local minima~\cite{horst2000introduction}. The two popular annealing methods include: (a) simulated annealing which is extensively used in machine learning algorithms~\cite{ackley1985learning,andrieu2003introduction}; and (b) quantum annealing which is at the core of emerging quantum computing processors like D-wave~\cite{dwave,devoret2013superconducting,kendon2010quantum,mohseni2017commercialize,pirandola2016unite}.  
Simulated annealing, as shown in FIG.~\ref{fig_globalopt}(a), searches for the lowest energy state (or global optimum of the target objective function) by gradually reducing the `temperature' of the system. As a result, during the initial stages of the procedure (high-temperature regime), the instantaneous solution (represented by the state of the electrons in FIG.~\ref{fig_globalopt}(a)) acquires sufficient energy to surmount the energy barriers surrounding the local minima \cite{kirkpatrick1983optimization}. As the temperature is cooled, the electrons become trapped in the global optimum (or lowest energy level) with a high-probability. However, due to its probabilistic nature, the simulated annealing approach does not provide any convergence guarantees and the heuristic choice of the annealing schedule strongly depends on the complexity of target objective function.  On the other hand, the quantum annealing approach~\cite{razavy2013quantum,johnson2011quantum}, as shown in FIG.~\ref{fig_globalopt}(b), uses a surrogate `driver' Hamiltonian function that has the identical lowest energy level as the original cost function, and evolves the latter over time adiabatically, to reach the global optimum, i.e., the state of minimum energy \cite{farhi2001quantum}. The energy-barriers separating the local energy minima are chosen to be sufficiently thin such that the instantaneous solution (or electrons in FIG.~\ref{fig_globalopt}(b)) can directly `tunnel through' to the global optimum from any locally optimal point. While quantum annealing algorithm overcomes some of the drawbacks of simulated annealing approaches and provides convergence guarantees, realizing and scaling such a system is cumbersome~\cite{barash1996low,dwave,harris2010experimental}, due to cryogenic operating conditions and due to the requirements on the thickness of the tunneling barriers.

In this paper we present an alternate approach for solving global optimization problems and is based on a novel dynamical system model that evolves
over a manifold defined by specific conservation constraints.  The basic principle is illustrated in FIG.~\ref{fig_globalopt}(c)  
where instead of directly optimizing the target objective function $q(\x))$, a `driver' function $h(\x,t)$ is allowed to evolve over an auxiliary Riemann manifold\cite{brockett1991dynamical} `$S$'  under the influence of $q(\x)$. At steady state $h(\x,t)$ converges to a Dirac-delta localized at the global minimum $\x^{*}$ ,i.e.,
\begin{equation}
h(\x,t)\overset{t\rightarrow\infty}{\longrightarrow} \delta(\x-\x^{*}),
\end{equation}
where $\delta(\x-\x^{*})$ is a Dirac-delta function at $\x\!=\!\x^{*}$, i.e., $\int \limits_{-\infty}^{+\infty}\delta(\x-\x^{*})d{\x}\!=\!1$. Note that during the process of
optimization, the solution (illustrated as electrons in FIG.~\ref{fig_globalopt}(c)) are free to move around on the auxiliary manifold `$S$'. 
\begin{algorithm*}[!htbp]
\caption{Main Result} 
\label{algo1}
\begin{algorithmic}
\vspace{0.2cm}
\STATE \underline{Objective:} 
\begin{equation}
 \underset{\x \in \mathbb{R}^M}{\text{minimize}} \quad q(\x),  \nonumber 
\end{equation}
where $q :\mathbb{R}^M \mapsto  \mathbb{R}$ is the target objective function has a single global minimum $\x^{*}$ such that $q(\x^{*})<q(\x), \forall \x \in D \subset \mathbb{R}^M$($q$ can have multiple local minima ) .\\
\vspace{0.1cm}
\STATE $\bullet$ Given a dynamical system described by the function $h(\x,t)$ which evolves according to:
\begin{align}
\tau \dfrac{\partial  h(\x,t)}{\partial t}+h(\x,t)=\dfrac{K(q(\x),\nu h(\x,t))}{\int _{\x} K(q(\x),\nu h(\x,t)) d{\x}} \nonumber
\end{align}
\vspace{0.1cm}
where $ K(.,.) $ is a functional with the form $K(q(\x),\nu h(\x,t))=h(\x,t)[\frac{1}{\nu}L\{q(\x),\nu h(\x,t)\} + \lambda]$, $L\{q(\x),\nu h(\x,t)\}$ is a monotonically increasing 
functional with respect to $h(\x,t)$, $\nu \in \mathbb{R}_+$, $\tau > 0$ is a system time-constant and $\lambda \in \mathcal{R}_+$ is a constant such that $K(.,.) > 0$.
\\
\vspace{0.2cm}

\STATE $\bullet$ Then, the following result holds:   
\begin{equation}
h(\x,t)\overset{t\rightarrow \infty}{\underset{\nu \to 0}{\longrightarrow}} \delta(\x-\x^{*}) \nonumber
\end{equation}
\end{algorithmic}
\end{algorithm*}
The description of a dynamical system that can achieve this evolution is summarized in Table~\ref{algo1} and the outline of the derivation of this main result is provided
at the end of this manuscript. The few important properties to note regarding the proposed optimization procedure: (a) The approach does not impose any differentiability constraints on the target objective function $q(\x)$ and hence 
in principle $q(\x)$ could be discrete valued; (b) The proposed method maps any multidimensional optimization problem to a one-dimensional space(since $h(\x,t)$ is scalar-valued), leading to a simplified treatment in a one-dimensional functional space; and (c) The proposed method
localizes the energy of the driver function at the global optimum, but does not directly provide the location of the Dirac-delta function. Therefore, akin to quantum
processors~\cite{mohseni2017commercialize,rabl2006hybrid}, a measurement process has to be used to infer the value of $\x^{*}$ using efficient random sampling techniques~\cite{cochran2007sampling,lo1999unconditional,shor2000simple}.
 
\begin{figure*}
\begin{minipage}[b]{0.48\textwidth}
\begin{center}
\includegraphics[scale=0.13,trim=4 4 4 4,clip]{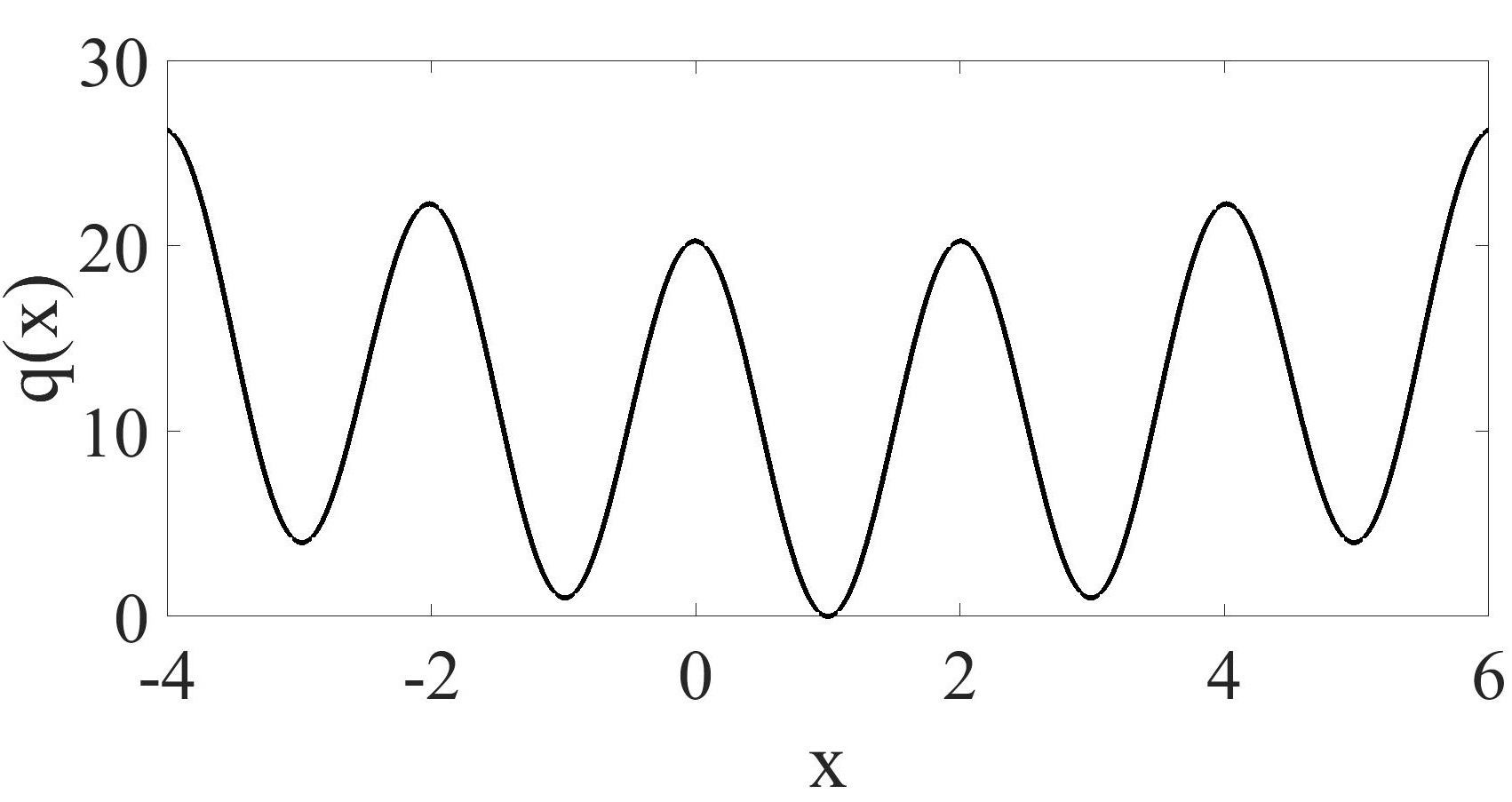}
\\
{\it (a)} \\
\end{center}
\end{minipage}
\begin{minipage}[b]{0.48\textwidth}
\begin{center}
\includegraphics[scale=0.13,trim=4 4 4 4,clip]{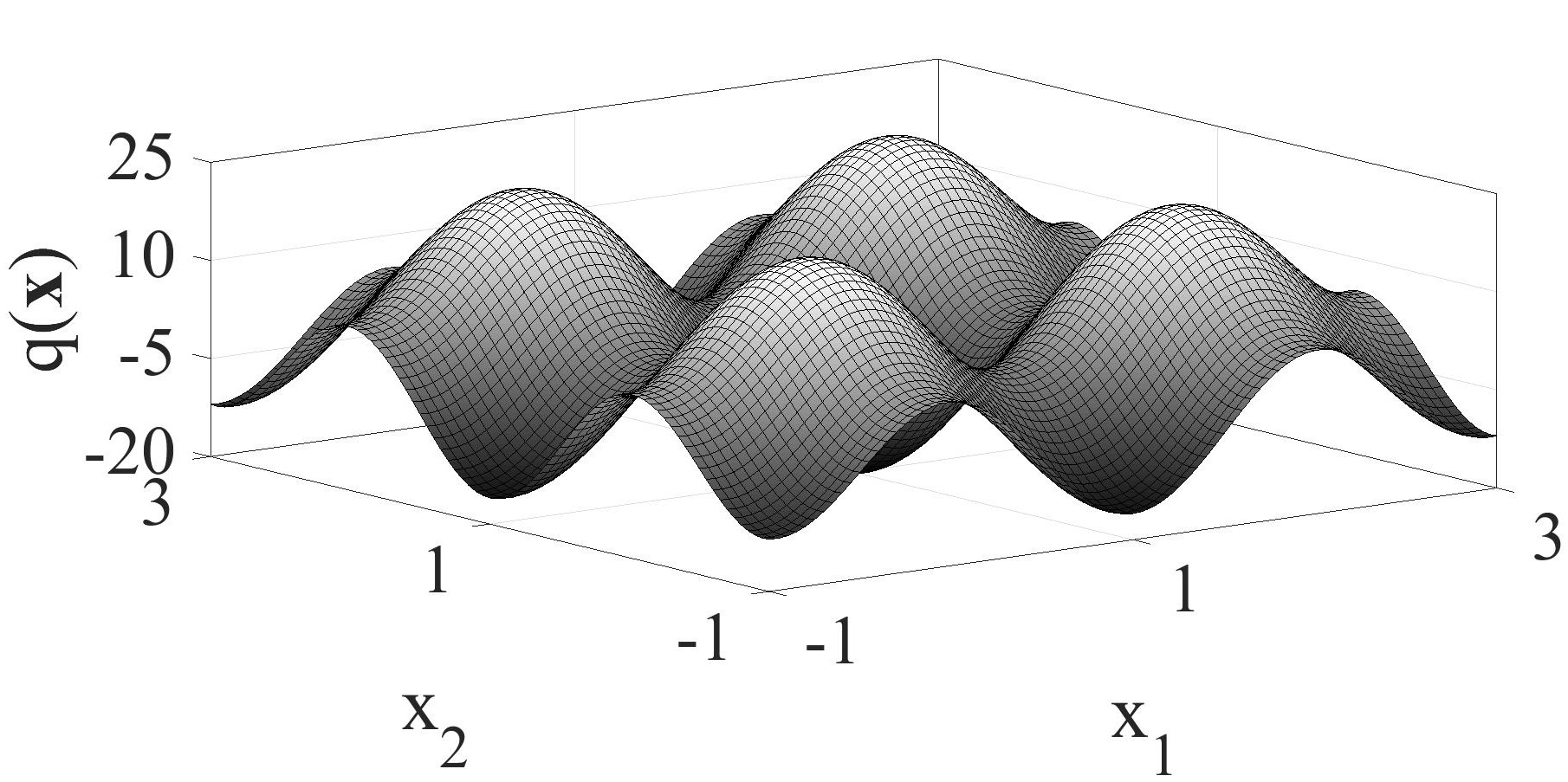}
\\
{\it (b)} \\
\end{center}
\end{minipage}
\begin{minipage}[b]{0.48\textwidth}
\begin{center}
\includegraphics[scale=0.13]{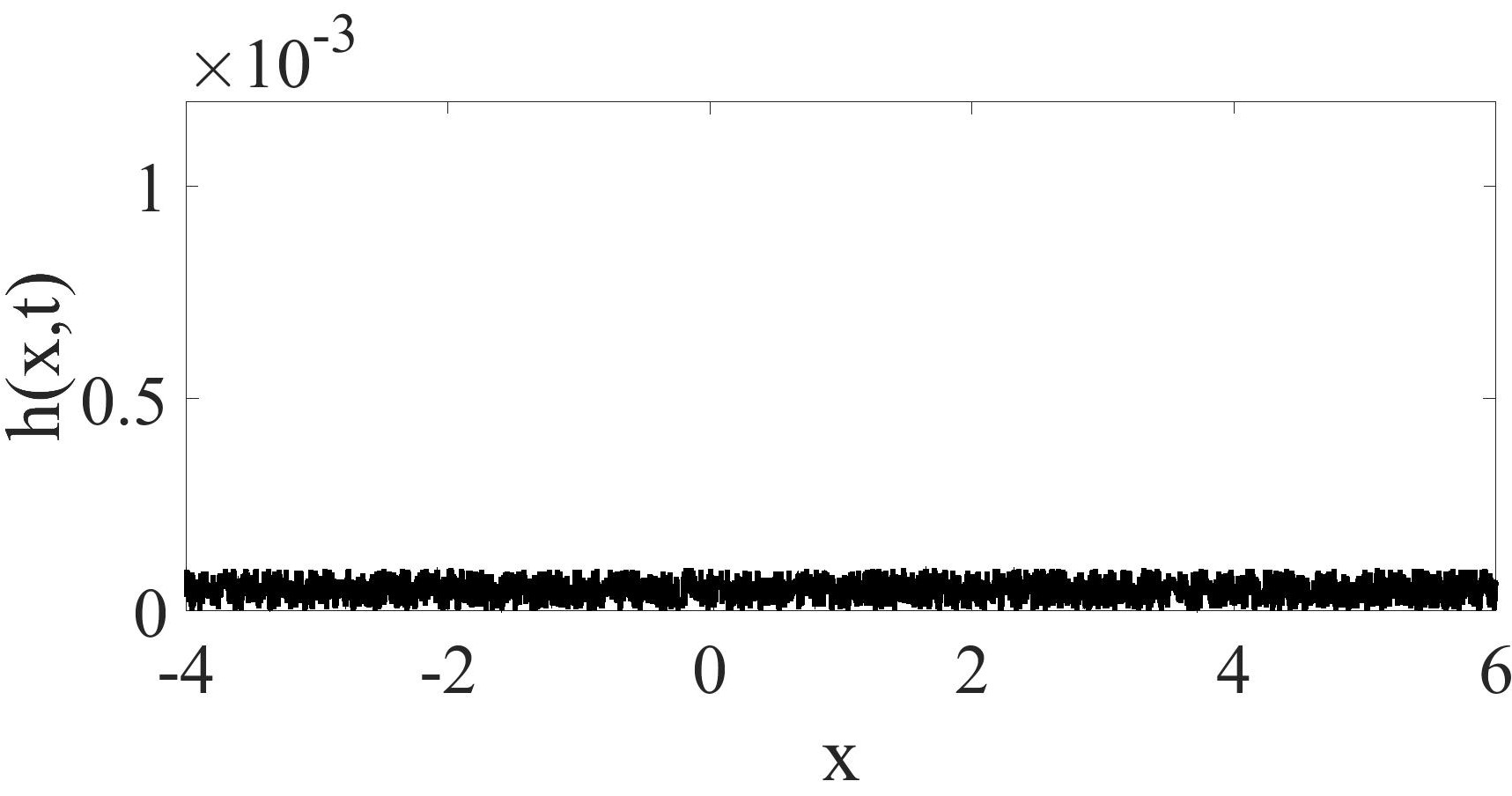} \\
{\it (c)} \\
\end{center}
\end{minipage}
\begin{minipage}[b]{0.48\textwidth}
\begin{center}
\includegraphics[scale=0.13]{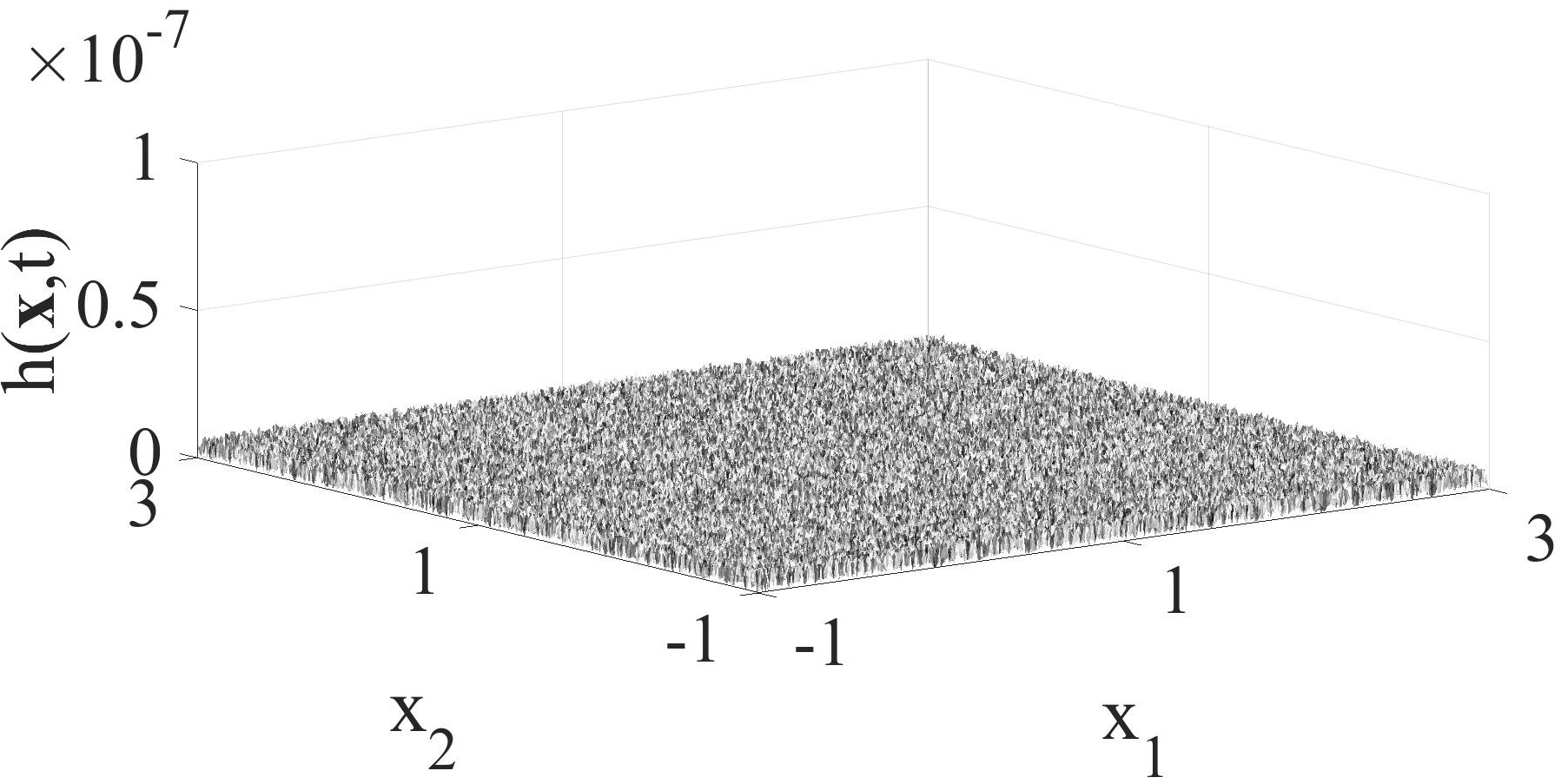} \\
{\it (d)} \\
\end{center}
\end{minipage}
\begin{minipage}[b]{0.48\textwidth}
\begin{center}
\includegraphics[scale=0.13]{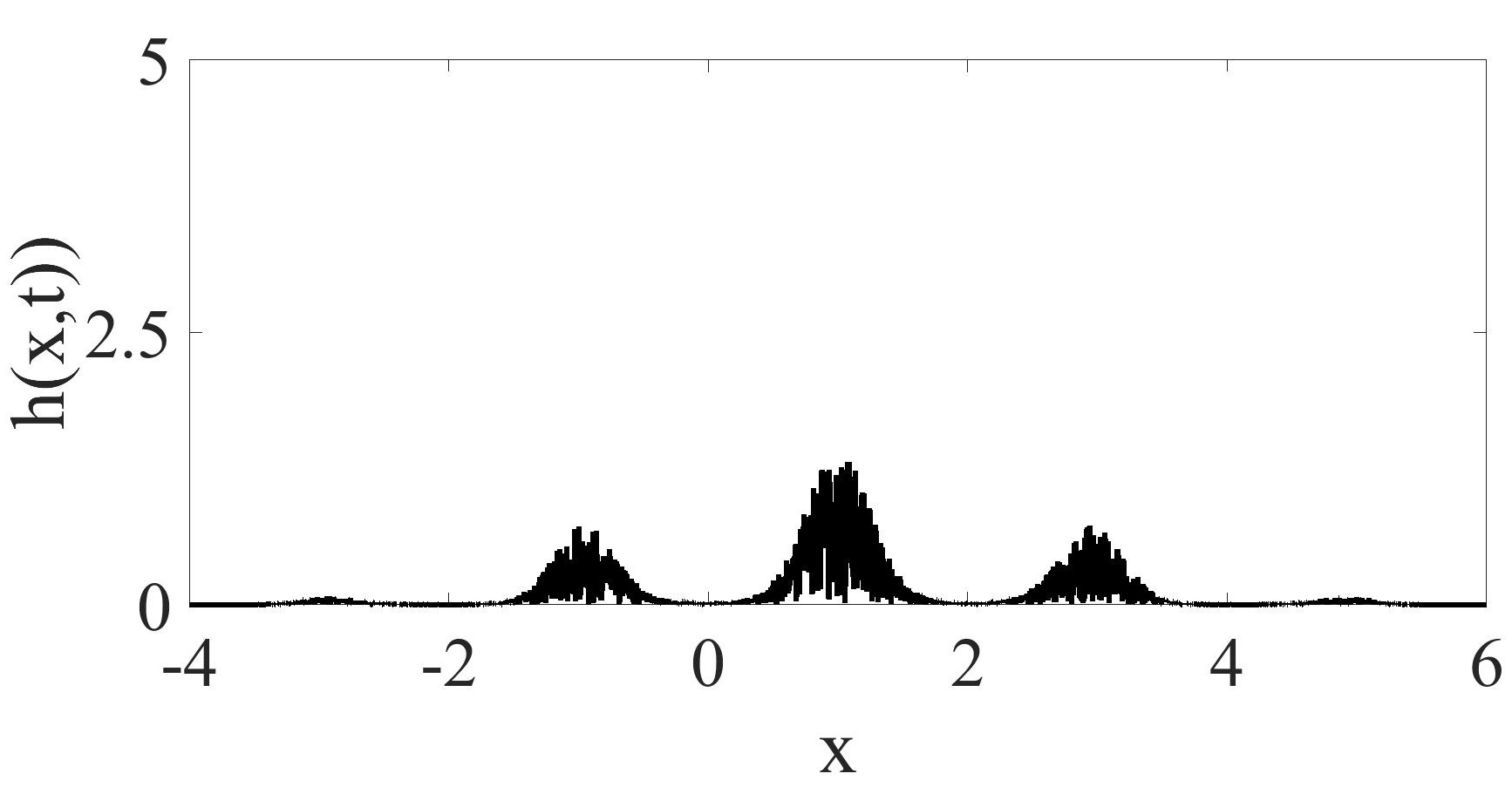} \\
{\it (e)} \\
\end{center}
\end{minipage}
\begin{minipage}[b]{0.48\textwidth}
\begin{center}
\includegraphics[scale=0.13]{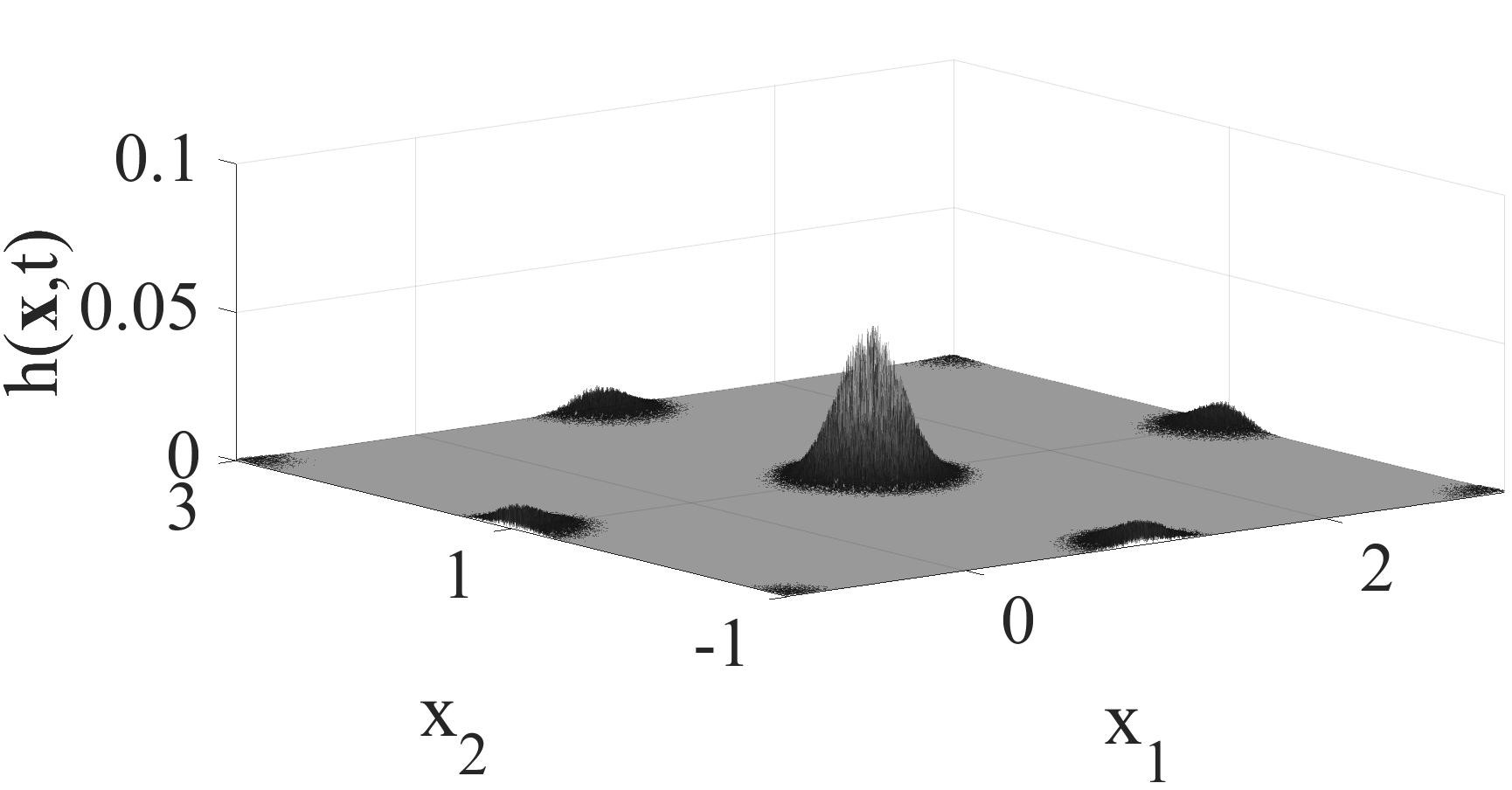} \\
{\it (f)} \\
\end{center}
\end{minipage}
\begin{minipage}[b]{0.48\textwidth}
\begin{center}
\includegraphics[scale=0.13]{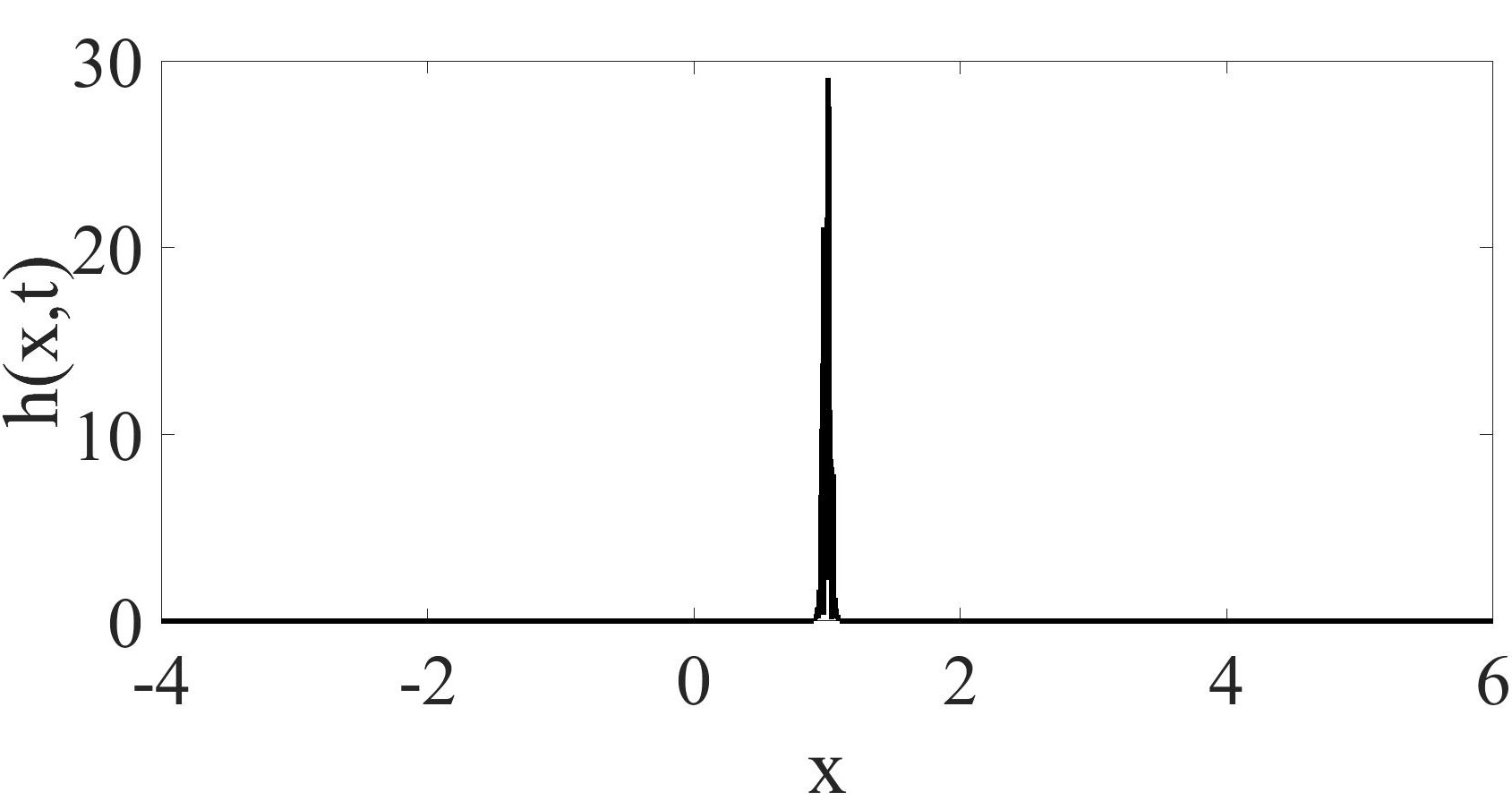} \\
{\it (g)} \\
\end{center}
\end{minipage}
\begin{minipage}[b]{0.48\textwidth}
\begin{center}
\includegraphics[scale=0.13]{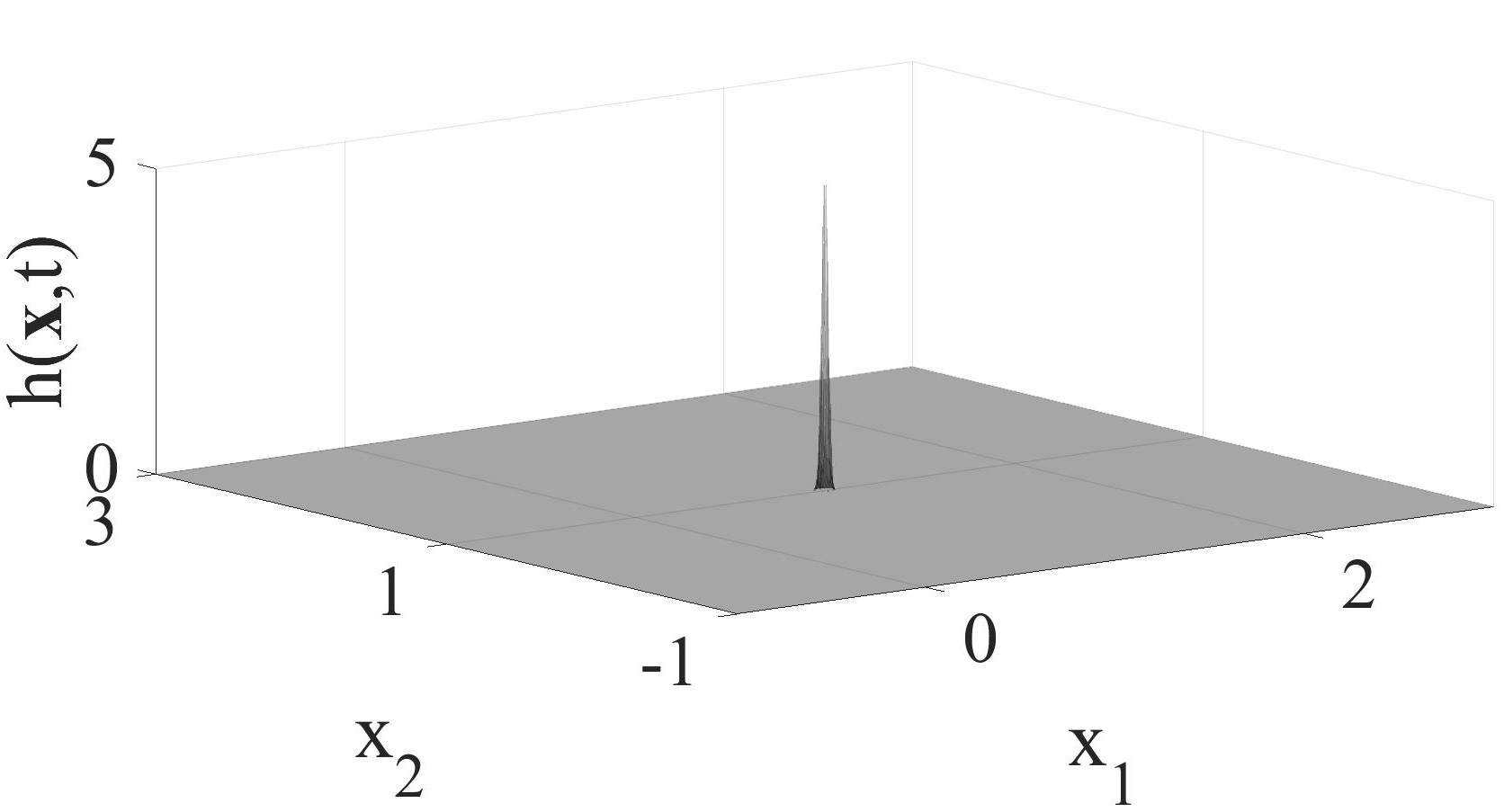} \\
{\it (h)} \\
\end{center}
\end{minipage}
\vspace{-0.5cm}
\caption{Illustration of the evolution of the distribution function $h(\x,t)$ when the dynamical system is used for optimizing (a) for a 1D rastrigin function and (b) for a 2D rastrigin function . In both cases, $h(\x,t)$ evolves from an initial random distribution (c), (d) to an intermediate (e),(f) to a Dirac-delta function located at $\x=\x^{*}$ in steady state (g),(h). }
\label{fig_evolve}
\end{figure*}

\begin{figure*}
\begin{center}
\includegraphics[page=1,scale=0.22,trim=4 4 4 4,clip]{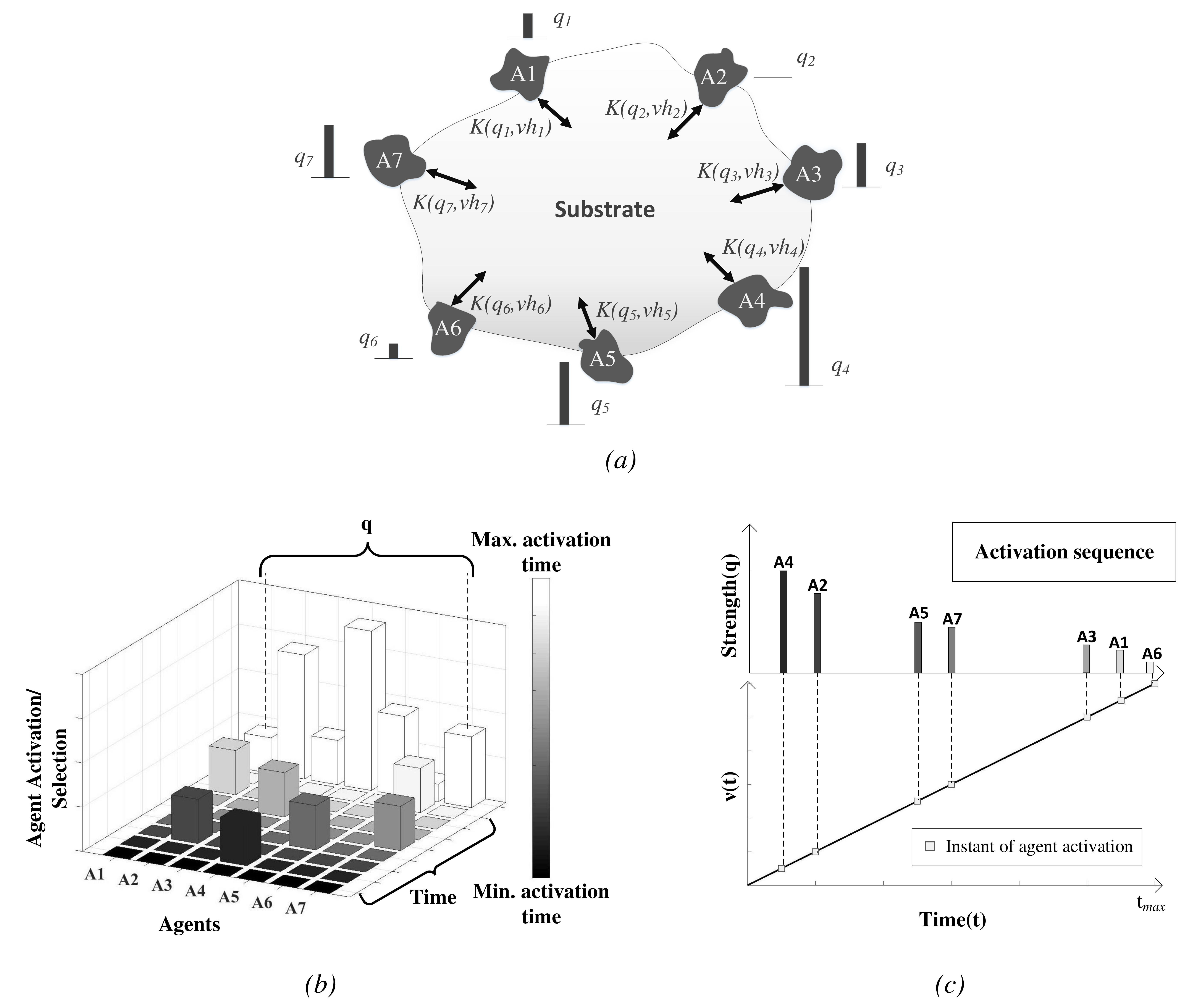}
\end{center}
\caption{Example of the dynamical system using for decentralized sorting: (a) Distributed agents, A1-A7 that only interact 
with its environment or substrate based on the local values $q_1$-$q_7$ and the interaction function $K(.,.)$ as described in Table~\ref{algo1}. The substrate simply integrates  the received $K(.,.)$ according to Table~\ref{algo1} which can then be accessed by each individual agents. (b) Linear-time sorting
algorithm where once the maximum is identified, the corresponding agent is deactivated from the optimization; and (c) Constant-time sorting algorithm
where $\nu$ is increased at a constant rate (slower than the time-constant $\tau$ in Table~\ref{algo1}) and each
of the agents are activated according to $q_1$-$q_7$. }
\label{fig_distributed}
\end{figure*}

\par We now apply the dynamical system model towards optimizing $M-$dimensional Rastrigin function which has served as a benchmark for verifying global optimization algorithms in literature\cite{potter1994cooperative}. A one-dimensional and two-dimensional variant of the Rastrigin function is shown in FIG. \ref{fig_evolve}(a) and (d), and given by $q(\x)= \sum \limits_{i=1}^M \{(x_i-1).^2-10\cos(\pi(x_i-1))\} + 10M$. Note that while $q(\x)$ has multiple local minima, 
it has only one global minimum at $\x^{*}=\mathbf{1}^M$. Also, while different forms of monotonic functional $L(.,.)$ could have been chosen, for this demonstration 
we have used $L\{q(\x),\nu h(\x,t)\}=[\nu h(\x,t) + q(\x))]$. The evolution of $h(\x,t)$ is shown in FIG.~\ref{fig_evolve} for three time instants and for both the one and two  dimensional examples. In both cases, starting from an initial random distribution(FIG. \ref{fig_evolve}(b),(f)), $h(\x,t)$ evolves over time (representative plots shown in FIG. \ref{fig_evolve}(c),(g)) to ultimately converge  to a Dirac-delta function centered at the global optimum $\x^{*}=\mathbf{1}^M$ in steady state (FIG. \ref{fig_evolve}(d),(h)). 
While the result in Table~\ref{algo1} has been reported for a continuous functional, the model could be discretized (with respect to $\x$) and applied
to decentralized optimization problems as well. In such a class of optimization problems, distributed entities (biological cells, particles or programs)
communicate only with their immediate environment (or substrate) through exchange of signals or chemicals~\cite{artavanis1999notch,lemmon2016dark,waters2005quorum,schlessinger2000cell}. Using the distributed exchange
of signals, the ensemble itself evolves to solve a global optimization task.  Examples of such decentralized optimization techniques include the nature-inspired heuristic methods like the evolutionary algorithms \cite{back1996evolutionary}, ant \cite{dorigo2006ant} and bee colony \cite{karaboga2007powerful} methods and particle swarm optimization \cite{kennedy2011particle}. Here we briefly demonstrate the proposed dynamical system model in Table~\ref{algo1} for decentralized
sorting. The driver functional $h(\x,t)$ is evaluated only at specific discrete points in $\x$ that correspond to the location of discrete agents $A1-A8$, as shown in FIG.~\ref{fig_distributed}. These
agents have access to localized values $q_1-q_8$ (which could represent concentration level of specific chemicals) and the objective is to
sort these values by exchanging local information with the substrate (shown in FIG.~\ref{fig_distributed}(a)). The substrate acts as a passive aggregator
of the local information that it receives from all the agents (or equivalently mixes the received signals). When the model in Table~\ref{algo1} is applied for
two specific cases, it leads to a linear-time complexity or a constant-time complexity sorting algorithm. Note that the time-complexity of 
standard sorting algorithms is $O(t \log t)$, with $t$ denoting absolute time.
\par {\it Practical issues and discussions:} 
The prescription of the dynamical system model in Table~\ref{algo1} is generic enough to be implementable in different types of devices and hardware.
The only aspect which could change between implementations would be the form of the functional $L(.,.)$. One possible practical realization
of this model could be using the physics of scanning transmission electron microscopes (STEM)\cite{binnig1982surface, williams2009transmission}. 
In this realization, the substrate `$S$' in FIG.~\ref{fig_globalopt}(c) could comprise of spatially delocalized mobile electrons, and the objective
function $q(\x)$ could be encoded using a conductive scaffold $C$ (connected to a higher potential than the substrate), as shown in FIG.~\ref{fig_globalopt}(c).
The constraint on the electrons on the substrate `$S$' could be imposed by throttling the flow of electrons (for instance using a constant current source).
As the scaffold and the substrate are brought closer to each other, the delocalized electrons on the substrate cluster towards the region $S^\prime$ that is closest to the minimum $C^\prime$, and at a close enough distance the electrons from
the substrate tunnel through the barrier between $S^\prime$ and $C^\prime$. By monitoring the time evolution of the electron density distribution over S(analogous to the `driver' function $h(x,t)$), the dynamical system can therefore be used
to solve the optimization problem. Note that because the
electrons are delocalized on the substrate S, the system is  always able to find the global minimum even if
the scaffold function C comprises of multiple minima, and the same argument holds true for the proposed model as well.
\vspace{-0.5cm}
\appendix*
\section{}
\vspace{-0.3cm}
Consider a functional optimization problem of the following form : 
\begin{gather}
\underset{f(\x): \mathbb{R}^{M} \mapsto \mathbb{R}_+} {\text{min}} H\{q(\x), f(\x)\}  \label{eq_lemma1}\\
\textit{s.t.}\quad \int \limits _{\x \in D \subset \mathbb{R}^{M}}f(\x)d{\x}=\nu,\quad f(\x)>0, \enskip \nu \in \mathbb{R}_+, \nonumber 
\end{gather}
where $q:\mathbb{R}^M \mapsto \mathbb{R}$ is any arbitrary function having multiple local minima, but a single global minimum at $\x^{*}$ (i.e., $q(\x^{*}) < q(\x), \enskip \forall \x \in D$), and $H:\mathbb{R}^M_+ \mapsto \mathbb{R}$ is a convex functional in $f(\x)$. Note that (\ref{eq_lemma1}) involves optimizing a convex cost functional $H$ over a convex domain, and thus has a unique solution.  By choosing $\nu \to 0$, we have $f(\x) \to 0 \enskip \forall \x \neq \x^{*}$. Thus, 
$f(\x)/\nu \to \delta(\x-\x^{*})$ such that $\int \limits_{-\infty}^{+\infty}\delta(\x-\x^{*})d{\x}\!=\!1$, i.e., the solution of (\ref{eq_lemma1}) is an impulse function at the minimal point of $q(\x)$. By denoting $h(\x) = f(\x)/\nu$, the functional optimization in~\ref{eq_lemma1} can be expressed as:
\begin{align}
\underset{\x \in\mathbb{R}^{M}} {\text{min}}q(\x) & \equiv \! \underset{h(\x): \mathbb{R}^{M} \mapsto \mathbb{R}_+} {\text{min}} H\{q(\x),\nu h(\x)\} \label{eq_lemma2}\\
& \textit{s.t.}\quad \int_{\x}h(\x)d{\x}=1,\quad h(\x)>0. \nonumber 
\end{align} 
\par We now show that the proposed dynamical system model in Table \ref{algo1} solves the optimization problem in~\ref{eq_lemma2}. For this
we will start with a discrete version of~\ref{eq_lemma2} and then extend it to the continuous domain. We will denote the discrete version of $H(q(\x),\nu h(\x))$ as $H(\p)$, with $p_i$ being the discretized version of $\nu h(\x)\Delta {\x}$, and the continuous domain $D_C=\{\int_{\x}h(\x)d{\x}=1,\quad h(\x)>0\}$ being transformed to the discrete domain $D=\{p_{i}: p_{i} \geq 0 \quad and \quad \sum\limits_{i=1}^{N} p_{i}=\nu\}$. Considering $H(\mathbf{p})$ to be  Lipschitz continuous on the compact domain $D$ with Lipschitz constant $ \lambda \! \in \! \mathbb{R}_+ $, i.e.,$ \lvert H(\mathbf{p_1})-H(\mathbf{p_2}) \rvert \le \lambda \lvert \lvert \mathbf{p_1}-\mathbf{p_2} \rvert \rvert_2, \forall \mathbf{p_1}, \mathbf{p_2} \in D$, it follows that $\lvert \dfrac{\partial H}{\partial p_i}\rvert < \lambda ,\quad \forall i=1,\ldots,N$. Now, if we define a function $G(\p)\!=\! -H(\p)+\lambda \sum \limits_{i=1}^N p_i$, then we have $\dfrac{\partial G}{\partial p_i}>0 ,\quad \forall i$, i.e., $G(\p)$ is monotonic, and the following thus holds:
\begin{align}
\underset{\p \in\mathbb{R}^{N}_+, \sum \limits _{i=1}^{N} p_i=\nu} {\text{min}} 
\! H(\p)=\underset{\p \in\mathbb{R}^{N}_+, \sum \limits _{i=1}^{N} p_i=\nu} {\text{max}} 
\!\{-H(\p)\} \nonumber\\
=\! \underset{\p \in\mathbb{R}^{N}_+, \sum \limits_{i=1}^{N} p_i=\nu} {\text{max}} \!G(\p)
\end{align}

\par Defining an auxiliary function $A(\p^t,\p)=\nabla G(\p^t)(\p-\p^t)$, we have $A(\p^t,\p^t)=0$. Additionally, if $A(\p^t,\p^{t+1})>0$ holds over 
$ D $, then it implies that $G(\p^{t+1})>G(\p^t)$ as well, since $G(\p)$ is monotonic. Thus, the function $G(\p)$ is bounded by a polynomial function $A(\p^t,\p)$, and can be optimized by locally optimizing $A(\p^t,\p)$ instead, as follows: 
\begin{gather}
 \underset{\p \in \mathbb{R}^N_+}{\text{maximize}} \quad A(\p,\p^t)  \label{eq_auxgrowth} \\
  \textit{s.t.} \quad \sum_{i=1}^{N} p_i=\nu ,\quad \p \succeq 0,\quad \nu \in \mathbb{R}_+  \nonumber .
\end{gather}
To do this we apply the celebrated Baum-Eagon growth transformation technique \cite{baum1968growth,gopalakrishnan1989generalization,svmgrowth}, where given $A:\mathbb{R}^{N}_+ \mapsto \mathbb{R}$ to be a polynomial with non-negative coefficients over $D$, if we define a continuous mapping $g:\mathbb{R}^{N}_+ \mapsto \mathbb{R}^{N}_+ $ as $g(p_{i})=\nu\dfrac{\tilde{g}_{i}(\mathbf{p}, \lambda)}{\tilde{g}(\mathbf{p},\lambda)} $,
where $\tilde g_{i}(\mathbf{p}, \lambda)=p_{i}(\frac{\partial A}{\partial p_{i}}+\lambda)$ and $\tilde g(\mathbf{p}, \lambda)=\sum\limits_{j=1}^{N} \tilde g_{j}(\mathbf{p}, \lambda)$, then $A(\p) \leq A(g(\p)), \quad \p \in D$.  $g$ thus maximizes $A$ over $D$, and is called a growth transformation \cite{gopalakrishnan1989generalization,svmgrowth} over $A$ on $D$. Additionally, for $0\leq\alpha\leq 1$, $A(\mathbf{p})\le A(\alpha\mathbf{p}+(1- \alpha)g(\mathbf{p}))$.
Thus, by minimizing the auxiliary function $A(\p)$, $H(\p)$ can thus be minimized using updates $g(p_i)=\nu \Bigg[\dfrac{p_i \Big(-\dfrac{\partial H}{\partial p_i}+\lambda\Big)}{\sum\limits_{j=1}^{N} p_j\Big(-\dfrac{\partial H}{\partial p_j}+\lambda \Big) }\Bigg], \quad \forall i=1,\ldots,N$.

\par Using the homotopically increasing property of growth transforms, i.e.,
\begin{equation}
p_i(t)=\alpha g(p_i(t-\Delta t))+(1- \alpha) p_i(t-\Delta t),\quad \forall i=1,\ldots,N
\nonumber
\end{equation}
and choosing $\alpha=\dfrac{1}{(\tau/\Delta t +1)}$, where $\tau$ can be thought of as the time constant of the ensemble system, we have 
\begin{equation}
\tau \bigg [\dfrac{\p(t)-\p(t-\Delta t)}{\Delta t} \bigg ]+\p(t)=g(\p(t-\Delta t)), \label{eq_dyn2}
\end{equation}
In the limiting case, when $\Delta t \to 0$, 
this reduces to the following continuous time dynamical system model: 
\begin{equation}
\tau \dfrac{d _ip(t)}{dt}+p_i(t)=g(p_i(t)),\quad \forall i=1, \ldots,N \label{eq_lemma5}
\end{equation}
where $g(p_i(t))=\nu \Bigg[ \dfrac{p_i(t)(-\dfrac{\partial H}{\partial p_i(t)}+\lambda)}{\sum \limits_j p_j(t)(-\dfrac{\partial H}{\partial p_j(t)}+\lambda)} \Bigg ]$.
 
\par  The continuous time framework (\ref{eq_lemma5}) can be extended to an infinite dimensional space \cite{rockafellar1974conjugate} and employed for optimization of functionals by taking $p_i\overset{N\rightarrow\infty}{\longrightarrow}\nu h(\x)\Delta \x ,\quad \x \in \mathbb{R}^M $, and considering the limiting case $\Delta \x \to 0$, to yield the following equation:
\begin{equation}
\tau \dfrac{\partial h(\x,t)}{\partial t}+h(\x,t)=g(h(\x,t)) ,\label{eq_dyncont}
\end{equation} 
where $g(h(\x,t))=\dfrac{h(\x,t) \Big\{-\dfrac{1}{\nu}\dfrac{\partial H}{\partial h(\x,t)}+\lambda\Big\}}{\int _{\x} h(\x,t) \Big\{-\dfrac{1}{\nu} \dfrac{\partial H}{\partial h(\x,t)}+\lambda \Big\} d{\x}}$.
The above dynamical system converges to the minimal point $\x^{*}$ of $H\{f(\x,t)\}$ in steady state, and can be expressed in a more compact and generic form as follows:
\begin{align}
\tau \dfrac{\partial  h(\x,t)}{\partial t}+h(\x,t)=\dfrac{K\{q(\x),\nu h(\x,t)\}}{\int _{\x} K\{q(\x),\nu h(\x,t)\} d{\x}}, 
\end{align} 
where $K(q(\x),\nu h(\x,t))=h(\x,t)[\frac{1}{\nu}L\{q(\x),\nu h(\x,t)\} + \lambda]$, and $L\{q(\x),\nu h(\x)\}$ is any monotonic functional in $h(\x)$.

\nocite{*}
\bibliography{references}

\begin{thebibliography}{39}%
\makeatletter
\providecommand \@ifxundefined [1]{%
 \@ifx{#1\undefined}
}%
\providecommand \@ifnum [1]{%
 \ifnum #1\expandafter \@firstoftwo
 \else \expandafter \@secondoftwo
 \fi
}%
\providecommand \@ifx [1]{%
 \ifx #1\expandafter \@firstoftwo
 \else \expandafter \@secondoftwo
 \fi
}%
\providecommand \natexlab [1]{#1}%
\providecommand \enquote  [1]{``#1''}%
\providecommand \bibnamefont  [1]{#1}%
\providecommand \bibfnamefont [1]{#1}%
\providecommand \citenamefont [1]{#1}%
\providecommand \href@noop [0]{\@secondoftwo}%
\providecommand \href [0]{\begingroup \@sanitize@url \@href}%
\providecommand \@href[1]{\@@startlink{#1}\@@href}%
\providecommand \@@href[1]{\endgroup#1\@@endlink}%
\providecommand \@sanitize@url [0]{\catcode `\\12\catcode `\$12\catcode
  `\&12\catcode `\#12\catcode `\^12\catcode `\_12\catcode `\%12\relax}%
\providecommand \@@startlink[1]{}%
\providecommand \@@endlink[0]{}%
\providecommand \url  [0]{\begingroup\@sanitize@url \@url }%
\providecommand \@url [1]{\endgroup\@href {#1}{\urlprefix }}%
\providecommand \urlprefix  [0]{URL }%
\providecommand \Eprint [0]{\href }%
\providecommand \doibase [0]{http://dx.doi.org/}%
\providecommand \selectlanguage [0]{\@gobble}%
\providecommand \bibinfo  [0]{\@secondoftwo}%
\providecommand \bibfield  [0]{\@secondoftwo}%
\providecommand \translation [1]{[#1]}%
\providecommand \BibitemOpen [0]{}%
\providecommand \bibitemStop [0]{}%
\providecommand \bibitemNoStop [0]{.\EOS\space}%
\providecommand \EOS [0]{\spacefactor3000\relax}%
\providecommand \BibitemShut  [1]{\csname bibitem#1\endcsname}%
\let\auto@bib@innerbib\@empty
\bibitem [{\citenamefont {Feynman}, \citenamefont {Leighton},\ and\
  \citenamefont {Sands}(2013)}]{feynman2013feynman}%
  \BibitemOpen
  \bibfield  {author} {\bibinfo {author} {\bibfnamefont {R.~P.}\ \bibnamefont
  {Feynman}}, \bibinfo {author} {\bibfnamefont {R.~B.}\ \bibnamefont
  {Leighton}}, \ and\ \bibinfo {author} {\bibfnamefont {M.}~\bibnamefont
  {Sands}},\ }\href@noop {} {}\ (\bibinfo  {publisher} {Basic books},\ \bibinfo
  {year} {2013})\BibitemShut {NoStop}%
\bibitem [{\citenamefont {Feynman}(1967)}]{feynman1967character}%
  \BibitemOpen
  \bibfield  {author} {\bibinfo {author} {\bibfnamefont {R.~P.}\ \bibnamefont
  {Feynman}},\ }\href@noop {} {}\ (\bibinfo  {publisher} {MIT press},\ \bibinfo
  {year} {1967})\BibitemShut {NoStop}%
\bibitem [{\citenamefont {Vergis}, \citenamefont {Steiglitz},\ and\
  \citenamefont {Dickinson}(1986)}]{vergis1986complexity}%
  \BibitemOpen
  \bibfield  {author} {\bibinfo {author} {\bibfnamefont {A.}~\bibnamefont
  {Vergis}}, \bibinfo {author} {\bibfnamefont {K.}~\bibnamefont {Steiglitz}}, \
  and\ \bibinfo {author} {\bibfnamefont {B.}~\bibnamefont {Dickinson}},\
  }\href@noop {} {\bibfield  {journal} {\bibinfo  {journal} {Mathematics and
  computers in simulation}\ }\textbf {\bibinfo {volume} {28}},\ \bibinfo
  {pages} {91} (\bibinfo {year} {1986})}\BibitemShut {NoStop}%
\bibitem [{\citenamefont {Hopfield}\ and\ \citenamefont
  {Tank}(1985)}]{hopfield1985neural}%
  \BibitemOpen
  \bibfield  {author} {\bibinfo {author} {\bibfnamefont {J.~J.}\ \bibnamefont
  {Hopfield}}\ and\ \bibinfo {author} {\bibfnamefont {D.~W.}\ \bibnamefont
  {Tank}},\ }\href@noop {} {\bibfield  {journal} {\bibinfo  {journal}
  {Biological cybernetics}\ }\textbf {\bibinfo {volume} {52}},\ \bibinfo
  {pages} {141} (\bibinfo {year} {1985})}\BibitemShut {NoStop}%
\bibitem [{\citenamefont {Hopfield}\ and\ \citenamefont
  {Herz}(1995)}]{hopfield1995rapid}%
  \BibitemOpen
  \bibfield  {author} {\bibinfo {author} {\bibfnamefont {J.~J.}\ \bibnamefont
  {Hopfield}}\ and\ \bibinfo {author} {\bibfnamefont {A.~V.}\ \bibnamefont
  {Herz}},\ }\href@noop {} {\bibfield  {journal} {\bibinfo  {journal}
  {Proceedings of the National Academy of Sciences}\ }\textbf {\bibinfo
  {volume} {92}},\ \bibinfo {pages} {6655} (\bibinfo {year}
  {1995})}\BibitemShut {NoStop}%
\bibitem [{\citenamefont {Horst}, \citenamefont {Pardalos},\ and\ \citenamefont
  {Van~Thoai}(2000)}]{horst2000introduction}%
  \BibitemOpen
  \bibfield  {author} {\bibinfo {author} {\bibfnamefont {R.}~\bibnamefont
  {Horst}}, \bibinfo {author} {\bibfnamefont {P.~M.}\ \bibnamefont {Pardalos}},
  \ and\ \bibinfo {author} {\bibfnamefont {N.}~\bibnamefont {Van~Thoai}},\
  }\href@noop {} {}\ (\bibinfo  {publisher} {Springer Science \& Business
  Media},\ \bibinfo {year} {2000})\BibitemShut {NoStop}%
\bibitem [{\citenamefont {Ackley}, \citenamefont {Hinton},\ and\ \citenamefont
  {Sejnowski}(1985)}]{ackley1985learning}%
  \BibitemOpen
  \bibfield  {author} {\bibinfo {author} {\bibfnamefont {D.~H.}\ \bibnamefont
  {Ackley}}, \bibinfo {author} {\bibfnamefont {G.~E.}\ \bibnamefont {Hinton}},
  \ and\ \bibinfo {author} {\bibfnamefont {T.~J.}\ \bibnamefont {Sejnowski}},\
  }\href@noop {} {\bibfield  {journal} {\bibinfo  {journal} {Cognitive
  science}\ }\textbf {\bibinfo {volume} {9}},\ \bibinfo {pages} {147} (\bibinfo
  {year} {1985})}\BibitemShut {NoStop}%
\bibitem [{\citenamefont {Andrieu}\ \emph {et~al.}(2003)\citenamefont
  {Andrieu}, \citenamefont {De~Freitas}, \citenamefont {Doucet},\ and\
  \citenamefont {Jordan}}]{andrieu2003introduction}%
  \BibitemOpen
  \bibfield  {author} {\bibinfo {author} {\bibfnamefont {C.}~\bibnamefont
  {Andrieu}}, \bibinfo {author} {\bibfnamefont {N.}~\bibnamefont {De~Freitas}},
  \bibinfo {author} {\bibfnamefont {A.}~\bibnamefont {Doucet}}, \ and\ \bibinfo
  {author} {\bibfnamefont {M.~I.}\ \bibnamefont {Jordan}},\ }\href@noop {}
  {\bibfield  {journal} {\bibinfo  {journal} {Machine learning}\ }\textbf
  {\bibinfo {volume} {50}},\ \bibinfo {pages} {5} (\bibinfo {year}
  {2003})}\BibitemShut {NoStop}%
\bibitem [{dwa(2017)}]{dwave}%
  \BibitemOpen
  \href@noop {} {}\bibinfo {type} {Tech. Rep.}\ (\bibinfo  {institution}
  {D-Wave Systems},\ \bibinfo {year} {2017})\BibitemShut {NoStop}%
\bibitem [{\citenamefont {Devoret}\ and\ \citenamefont
  {Schoelkopf}(2013)}]{devoret2013superconducting}%
  \BibitemOpen
  \bibfield  {author} {\bibinfo {author} {\bibfnamefont {M.~H.}\ \bibnamefont
  {Devoret}}\ and\ \bibinfo {author} {\bibfnamefont {R.~J.}\ \bibnamefont
  {Schoelkopf}},\ }\href@noop {} {\bibfield  {journal} {\bibinfo  {journal}
  {Science}\ }\textbf {\bibinfo {volume} {339}},\ \bibinfo {pages} {1169}
  (\bibinfo {year} {2013})}\BibitemShut {NoStop}%
\bibitem [{\citenamefont {Kendon}, \citenamefont {Nemoto},\ and\ \citenamefont
  {Munro}(2010)}]{kendon2010quantum}%
  \BibitemOpen
  \bibfield  {author} {\bibinfo {author} {\bibfnamefont {V.~M.}\ \bibnamefont
  {Kendon}}, \bibinfo {author} {\bibfnamefont {K.}~\bibnamefont {Nemoto}}, \
  and\ \bibinfo {author} {\bibfnamefont {W.~J.}\ \bibnamefont {Munro}},\
  }\href@noop {} {\bibfield  {journal} {\bibinfo  {journal} {Philosophical
  Transactions of the Royal Society of London A: Mathematical, Physical and
  Engineering Sciences}\ }\textbf {\bibinfo {volume} {368}},\ \bibinfo {pages}
  {3609} (\bibinfo {year} {2010})}\BibitemShut {NoStop}%
\bibitem [{\citenamefont {Mohseni}\ \emph {et~al.}(2017)\citenamefont
  {Mohseni}, \citenamefont {Read}, \citenamefont {Neven}, \citenamefont
  {Boixo}, \citenamefont {Denchev}, \citenamefont {Babbush}, \citenamefont
  {Fowler}, \citenamefont {Smelyanskiy}, \citenamefont {Martinis} \emph
  {et~al.}}]{mohseni2017commercialize}%
  \BibitemOpen
  \bibfield  {author} {\bibinfo {author} {\bibfnamefont {M.}~\bibnamefont
  {Mohseni}}, \bibinfo {author} {\bibfnamefont {P.}~\bibnamefont {Read}},
  \bibinfo {author} {\bibfnamefont {H.}~\bibnamefont {Neven}}, \bibinfo
  {author} {\bibfnamefont {S.}~\bibnamefont {Boixo}}, \bibinfo {author}
  {\bibfnamefont {V.}~\bibnamefont {Denchev}}, \bibinfo {author} {\bibfnamefont
  {R.}~\bibnamefont {Babbush}}, \bibinfo {author} {\bibfnamefont
  {A.}~\bibnamefont {Fowler}}, \bibinfo {author} {\bibfnamefont
  {V.}~\bibnamefont {Smelyanskiy}}, \bibinfo {author} {\bibfnamefont
  {J.}~\bibnamefont {Martinis}},  \emph {et~al.},\ }\href@noop {} {} (\bibinfo
  {year} {2017})\BibitemShut {NoStop}%
\bibitem [{\citenamefont {Pirandola}\ and\ \citenamefont
  {Braunstein}(2016)}]{pirandola2016unite}%
  \BibitemOpen
  \bibfield  {author} {\bibinfo {author} {\bibfnamefont {S.}~\bibnamefont
  {Pirandola}}\ and\ \bibinfo {author} {\bibfnamefont {S.~L.}\ \bibnamefont
  {Braunstein}},\ }\href@noop {} {\bibfield  {journal} {\bibinfo  {journal}
  {Nature}\ ,\ \bibinfo {pages} {169}} (\bibinfo {year} {2016})}\BibitemShut
  {NoStop}%
\bibitem [{\citenamefont {Kirkpatrick}\ \emph {et~al.}(1983)\citenamefont
  {Kirkpatrick}, \citenamefont {Gelatt}, \citenamefont {Vecchi} \emph
  {et~al.}}]{kirkpatrick1983optimization}%
  \BibitemOpen
  \bibfield  {author} {\bibinfo {author} {\bibfnamefont {S.}~\bibnamefont
  {Kirkpatrick}}, \bibinfo {author} {\bibfnamefont {C.~D.}\ \bibnamefont
  {Gelatt}}, \bibinfo {author} {\bibfnamefont {M.~P.}\ \bibnamefont {Vecchi}},
  \emph {et~al.},\ }\href@noop {} {\bibfield  {journal} {\bibinfo  {journal}
  {Science}\ }\textbf {\bibinfo {volume} {220}},\ \bibinfo {pages} {671}
  (\bibinfo {year} {1983})}\BibitemShut {NoStop}%
\bibitem [{\citenamefont {Razavy}(2013)}]{razavy2013quantum}%
  \BibitemOpen
  \bibfield  {author} {\bibinfo {author} {\bibfnamefont {M.}~\bibnamefont
  {Razavy}},\ }\href@noop {} {\emph {\bibinfo {title} {Quantum theory of
  tunneling}}}\ (\bibinfo  {publisher} {World Scientific},\ \bibinfo {year}
  {2013})\BibitemShut {NoStop}%
\bibitem [{\citenamefont {Johnson}\ \emph {et~al.}(2011)\citenamefont
  {Johnson}, \citenamefont {Amin}, \citenamefont {Gildert}, \citenamefont
  {Lanting}, \citenamefont {Hamze}, \citenamefont {Dickson}, \citenamefont
  {Harris}, \citenamefont {Berkley}, \citenamefont {Johansson}, \citenamefont
  {Bunyk} \emph {et~al.}}]{johnson2011quantum}%
  \BibitemOpen
  \bibfield  {author} {\bibinfo {author} {\bibfnamefont {M.~W.}\ \bibnamefont
  {Johnson}}, \bibinfo {author} {\bibfnamefont {M.~H.}\ \bibnamefont {Amin}},
  \bibinfo {author} {\bibfnamefont {S.}~\bibnamefont {Gildert}}, \bibinfo
  {author} {\bibfnamefont {T.}~\bibnamefont {Lanting}}, \bibinfo {author}
  {\bibfnamefont {F.}~\bibnamefont {Hamze}}, \bibinfo {author} {\bibfnamefont
  {N.}~\bibnamefont {Dickson}}, \bibinfo {author} {\bibfnamefont
  {R.}~\bibnamefont {Harris}}, \bibinfo {author} {\bibfnamefont {A.~J.}\
  \bibnamefont {Berkley}}, \bibinfo {author} {\bibfnamefont {J.}~\bibnamefont
  {Johansson}}, \bibinfo {author} {\bibfnamefont {P.}~\bibnamefont {Bunyk}},
  \emph {et~al.},\ }\href@noop {} {\bibfield  {journal} {\bibinfo  {journal}
  {Nature}\ }\textbf {\bibinfo {volume} {473}} (\bibinfo {year}
  {2011})}\BibitemShut {NoStop}%
\bibitem [{\citenamefont {Farhi}\ \emph {et~al.}(2001)\citenamefont {Farhi},
  \citenamefont {Goldstone}, \citenamefont {Gutmann}, \citenamefont {Lapan},
  \citenamefont {Lundgren},\ and\ \citenamefont {Preda}}]{farhi2001quantum}%
  \BibitemOpen
  \bibfield  {author} {\bibinfo {author} {\bibfnamefont {E.}~\bibnamefont
  {Farhi}}, \bibinfo {author} {\bibfnamefont {J.}~\bibnamefont {Goldstone}},
  \bibinfo {author} {\bibfnamefont {S.}~\bibnamefont {Gutmann}}, \bibinfo
  {author} {\bibfnamefont {J.}~\bibnamefont {Lapan}}, \bibinfo {author}
  {\bibfnamefont {A.}~\bibnamefont {Lundgren}}, \ and\ \bibinfo {author}
  {\bibfnamefont {D.}~\bibnamefont {Preda}},\ }\href@noop {} {\bibfield
  {journal} {\bibinfo  {journal} {Science}\ }\textbf {\bibinfo {volume}
  {292}},\ \bibinfo {pages} {472} (\bibinfo {year} {2001})}\BibitemShut
  {NoStop}%
\bibitem [{\citenamefont {Barash}, \citenamefont {Burkhardt},\ and\
  \citenamefont {Rainer}(1996)}]{barash1996low}%
  \BibitemOpen
  \bibfield  {author} {\bibinfo {author} {\bibfnamefont {Y.~S.}\ \bibnamefont
  {Barash}}, \bibinfo {author} {\bibfnamefont {H.}~\bibnamefont {Burkhardt}}, \
  and\ \bibinfo {author} {\bibfnamefont {D.}~\bibnamefont {Rainer}},\
  }\href@noop {} {\bibfield  {journal} {\bibinfo  {journal} {Physical review
  letters}\ }\textbf {\bibinfo {volume} {77}},\ \bibinfo {pages} {4070}
  (\bibinfo {year} {1996})}\BibitemShut {NoStop}%
\bibitem [{\citenamefont {Harris}\ \emph {et~al.}(2010)\citenamefont {Harris},
  \citenamefont {Johnson}, \citenamefont {Lanting}, \citenamefont {Berkley},
  \citenamefont {Johansson}, \citenamefont {Bunyk}, \citenamefont {Tolkacheva},
  \citenamefont {Ladizinsky}, \citenamefont {Ladizinsky}, \citenamefont {Oh}
  \emph {et~al.}}]{harris2010experimental}%
  \BibitemOpen
  \bibfield  {author} {\bibinfo {author} {\bibfnamefont {R.}~\bibnamefont
  {Harris}}, \bibinfo {author} {\bibfnamefont {M.}~\bibnamefont {Johnson}},
  \bibinfo {author} {\bibfnamefont {T.}~\bibnamefont {Lanting}}, \bibinfo
  {author} {\bibfnamefont {A.}~\bibnamefont {Berkley}}, \bibinfo {author}
  {\bibfnamefont {J.}~\bibnamefont {Johansson}}, \bibinfo {author}
  {\bibfnamefont {P.}~\bibnamefont {Bunyk}}, \bibinfo {author} {\bibfnamefont
  {E.}~\bibnamefont {Tolkacheva}}, \bibinfo {author} {\bibfnamefont
  {E.}~\bibnamefont {Ladizinsky}}, \bibinfo {author} {\bibfnamefont
  {N.}~\bibnamefont {Ladizinsky}}, \bibinfo {author} {\bibfnamefont
  {T.}~\bibnamefont {Oh}},  \emph {et~al.},\ }\href@noop {} {\bibfield
  {journal} {\bibinfo  {journal} {Physical Review B}\ }\textbf {\bibinfo
  {volume} {024511}} (\bibinfo {year} {2010})}\BibitemShut {NoStop}%
\bibitem [{\citenamefont {Brockett}(1991)}]{brockett1991dynamical}%
  \BibitemOpen
  \bibfield  {author} {\bibinfo {author} {\bibfnamefont {R.~W.}\ \bibnamefont
  {Brockett}},\ }\href@noop {} {\bibfield  {journal} {\bibinfo  {journal}
  {Linear Algebra and its applications}\ }\textbf {\bibinfo {volume} {146}},\
  \bibinfo {pages} {79} (\bibinfo {year} {1991})}\BibitemShut {NoStop}%
\bibitem [{\citenamefont {Rabl}\ \emph {et~al.}(2006)\citenamefont {Rabl},
  \citenamefont {DeMille}, \citenamefont {Doyle}, \citenamefont {Lukin},
  \citenamefont {Schoelkopf},\ and\ \citenamefont {Zoller}}]{rabl2006hybrid}%
  \BibitemOpen
  \bibfield  {author} {\bibinfo {author} {\bibfnamefont {P.}~\bibnamefont
  {Rabl}}, \bibinfo {author} {\bibfnamefont {D.}~\bibnamefont {DeMille}},
  \bibinfo {author} {\bibfnamefont {J.~M.}\ \bibnamefont {Doyle}}, \bibinfo
  {author} {\bibfnamefont {M.~D.}\ \bibnamefont {Lukin}}, \bibinfo {author}
  {\bibfnamefont {R.}~\bibnamefont {Schoelkopf}}, \ and\ \bibinfo {author}
  {\bibfnamefont {P.}~\bibnamefont {Zoller}},\ }\href@noop {} {\bibfield
  {journal} {\bibinfo  {journal} {Physical review letters}\ }\textbf {\bibinfo
  {volume} {97}},\ \bibinfo {pages} {033003} (\bibinfo {year}
  {2006})}\BibitemShut {NoStop}%
\bibitem [{\citenamefont {Cochran}(2007)}]{cochran2007sampling}%
  \BibitemOpen
  \bibfield  {author} {\bibinfo {author} {\bibfnamefont {W.~G.}\ \bibnamefont
  {Cochran}},\ }\href@noop {} {}\ (\bibinfo  {publisher} {John Wiley \& Sons},\
  \bibinfo {year} {2007})\BibitemShut {NoStop}%
\bibitem [{\citenamefont {Lo}\ and\ \citenamefont
  {Chau}(1999)}]{lo1999unconditional}%
  \BibitemOpen
  \bibfield  {author} {\bibinfo {author} {\bibfnamefont {H.-K.}\ \bibnamefont
  {Lo}}\ and\ \bibinfo {author} {\bibfnamefont {H.~F.}\ \bibnamefont {Chau}},\
  }\href@noop {} {\bibfield  {journal} {\bibinfo  {journal} {Science}\ }\textbf
  {\bibinfo {volume} {283}},\ \bibinfo {pages} {2050} (\bibinfo {year}
  {1999})}\BibitemShut {NoStop}%
\bibitem [{\citenamefont {Shor}\ and\ \citenamefont
  {Preskill}(2000)}]{shor2000simple}%
  \BibitemOpen
  \bibfield  {author} {\bibinfo {author} {\bibfnamefont {P.~W.}\ \bibnamefont
  {Shor}}\ and\ \bibinfo {author} {\bibfnamefont {J.}~\bibnamefont
  {Preskill}},\ }\href@noop {} {\bibfield  {journal} {\bibinfo  {journal}
  {Physical review letters}\ }\textbf {\bibinfo {volume} {85}},\ \bibinfo
  {pages} {441} (\bibinfo {year} {2000})}\BibitemShut {NoStop}%
\bibitem [{\citenamefont {Potter}\ and\ \citenamefont
  {De~Jong}(1994)}]{potter1994cooperative}%
  \BibitemOpen
  \bibfield  {author} {\bibinfo {author} {\bibfnamefont {M.~A.}\ \bibnamefont
  {Potter}}\ and\ \bibinfo {author} {\bibfnamefont {K.~A.}\ \bibnamefont
  {De~Jong}},\ }in\ \href@noop {} {\emph {\bibinfo {booktitle} {International
  Conference on Parallel Problem Solving from Nature}}}\ (\bibinfo
  {organization} {Springer},\ \bibinfo {year} {1994})\ p.\ \bibinfo {pages}
  {249}\BibitemShut {NoStop}%
\bibitem [{\citenamefont {Artavanis-Tsakonas}, \citenamefont {Rand},\ and\
  \citenamefont {Lake}(1999)}]{artavanis1999notch}%
  \BibitemOpen
  \bibfield  {author} {\bibinfo {author} {\bibfnamefont {S.}~\bibnamefont
  {Artavanis-Tsakonas}}, \bibinfo {author} {\bibfnamefont {M.~D.}\ \bibnamefont
  {Rand}}, \ and\ \bibinfo {author} {\bibfnamefont {R.~J.}\ \bibnamefont
  {Lake}},\ }\href@noop {} {\bibfield  {journal} {\bibinfo  {journal}
  {Science}\ }\textbf {\bibinfo {volume} {284}},\ \bibinfo {pages} {770}
  (\bibinfo {year} {1999})}\BibitemShut {NoStop}%
\bibitem [{\citenamefont {Lemmon}\ \emph {et~al.}(2016)\citenamefont {Lemmon},
  \citenamefont {Freed}, \citenamefont {Schlessinger},\ and\ \citenamefont
  {Kiyatkin}}]{lemmon2016dark}%
  \BibitemOpen
  \bibfield  {author} {\bibinfo {author} {\bibfnamefont {M.~A.}\ \bibnamefont
  {Lemmon}}, \bibinfo {author} {\bibfnamefont {D.~M.}\ \bibnamefont {Freed}},
  \bibinfo {author} {\bibfnamefont {J.}~\bibnamefont {Schlessinger}}, \ and\
  \bibinfo {author} {\bibfnamefont {A.}~\bibnamefont {Kiyatkin}},\ }\href@noop
  {} {\bibfield  {journal} {\bibinfo  {journal} {Cell}\ }\textbf {\bibinfo
  {volume} {164}},\ \bibinfo {pages} {1172} (\bibinfo {year}
  {2016})}\BibitemShut {NoStop}%
\bibitem [{\citenamefont {Waters}\ and\ \citenamefont
  {Bassler}(2005)}]{waters2005quorum}%
  \BibitemOpen
  \bibfield  {author} {\bibinfo {author} {\bibfnamefont {C.~M.}\ \bibnamefont
  {Waters}}\ and\ \bibinfo {author} {\bibfnamefont {B.~L.}\ \bibnamefont
  {Bassler}},\ }\href@noop {} {\bibfield  {journal} {\bibinfo  {journal} {Annu.
  Rev. Cell Dev. Biol.}\ }\textbf {\bibinfo {volume} {21}},\ \bibinfo {pages}
  {319} (\bibinfo {year} {2005})}\BibitemShut {NoStop}%
\bibitem [{\citenamefont {Schlessinger}(2000)}]{schlessinger2000cell}%
  \BibitemOpen
  \bibfield  {author} {\bibinfo {author} {\bibfnamefont {J.}~\bibnamefont
  {Schlessinger}},\ }\href@noop {} {\bibfield  {journal} {\bibinfo  {journal}
  {Cell}\ }\textbf {\bibinfo {volume} {103}},\ \bibinfo {pages} {211} (\bibinfo
  {year} {2000})}\BibitemShut {NoStop}%
\bibitem [{\citenamefont {Back}(1996)}]{back1996evolutionary}%
  \BibitemOpen
  \bibfield  {author} {\bibinfo {author} {\bibfnamefont {T.}~\bibnamefont
  {Back}},\ }\href@noop {} {}\ (\bibinfo  {publisher} {Oxford university
  press},\ \bibinfo {year} {1996})\BibitemShut {NoStop}%
\bibitem [{\citenamefont {Dorigo}, \citenamefont {Birattari},\ and\
  \citenamefont {Stutzle}(2006)}]{dorigo2006ant}%
  \BibitemOpen
  \bibfield  {author} {\bibinfo {author} {\bibfnamefont {M.}~\bibnamefont
  {Dorigo}}, \bibinfo {author} {\bibfnamefont {M.}~\bibnamefont {Birattari}}, \
  and\ \bibinfo {author} {\bibfnamefont {T.}~\bibnamefont {Stutzle}},\
  }\href@noop {} {\bibfield  {journal} {\bibinfo  {journal} {IEEE computational
  intelligence magazine}\ }\textbf {\bibinfo {volume} {1}},\ \bibinfo {pages}
  {28} (\bibinfo {year} {2006})}\BibitemShut {NoStop}%
\bibitem [{\citenamefont {Karaboga}\ and\ \citenamefont
  {Basturk}(2007)}]{karaboga2007powerful}%
  \BibitemOpen
  \bibfield  {author} {\bibinfo {author} {\bibfnamefont {D.}~\bibnamefont
  {Karaboga}}\ and\ \bibinfo {author} {\bibfnamefont {B.}~\bibnamefont
  {Basturk}},\ }\href@noop {} {\bibfield  {journal} {\bibinfo  {journal}
  {Journal of global optimization}\ }\textbf {\bibinfo {volume} {39}},\
  \bibinfo {pages} {459} (\bibinfo {year} {2007})}\BibitemShut {NoStop}%
\bibitem [{\citenamefont {Kennedy}(2011)}]{kennedy2011particle}%
  \BibitemOpen
  \bibfield  {author} {\bibinfo {author} {\bibfnamefont {J.}~\bibnamefont
  {Kennedy}},\ }in\ \href@noop {} {\emph {\bibinfo {booktitle} {Encyclopedia of
  machine learning}}}\ (\bibinfo  {publisher} {Springer},\ \bibinfo {year}
  {2011})\ p.\ \bibinfo {pages} {760}\BibitemShut {NoStop}%
\bibitem [{\citenamefont {Binnig}\ \emph {et~al.}(1982)\citenamefont {Binnig},
  \citenamefont {Rohrer}, \citenamefont {Gerber},\ and\ \citenamefont
  {Weibel}}]{binnig1982surface}%
  \BibitemOpen
  \bibfield  {author} {\bibinfo {author} {\bibfnamefont {G.}~\bibnamefont
  {Binnig}}, \bibinfo {author} {\bibfnamefont {H.}~\bibnamefont {Rohrer}},
  \bibinfo {author} {\bibfnamefont {C.}~\bibnamefont {Gerber}}, \ and\ \bibinfo
  {author} {\bibfnamefont {E.}~\bibnamefont {Weibel}},\ }\href@noop {}
  {\bibfield  {journal} {\bibinfo  {journal} {Physical review letters}\
  }\textbf {\bibinfo {volume} {49}},\ \bibinfo {pages} {57} (\bibinfo {year}
  {1982})}\BibitemShut {NoStop}%
\bibitem [{\citenamefont {Williams}\ and\ \citenamefont
  {Carter}(2009)}]{williams2009transmission}%
  \BibitemOpen
  \bibfield  {author} {\bibinfo {author} {\bibfnamefont {D.~B.}\ \bibnamefont
  {Williams}}\ and\ \bibinfo {author} {\bibfnamefont {C.~B.}\ \bibnamefont
  {Carter}},\ }\href@noop {} {\bibfield  {journal} {\bibinfo  {journal}
  {Transmission electron microscopy}\ }\textbf {\bibinfo {volume} {3}}
  (\bibinfo {year} {2009})}\BibitemShut {NoStop}%
\bibitem [{\citenamefont {Baum}\ and\ \citenamefont
  {Sell}(1968)}]{baum1968growth}%
  \BibitemOpen
  \bibfield  {author} {\bibinfo {author} {\bibfnamefont {L.~E.}\ \bibnamefont
  {Baum}}\ and\ \bibinfo {author} {\bibfnamefont {G.~R.}\ \bibnamefont
  {Sell}},\ }\href@noop {} {\bibfield  {journal} {\bibinfo  {journal} {Pacific
  J. Math}\ }\textbf {\bibinfo {volume} {27}},\ \bibinfo {pages} {211}
  (\bibinfo {year} {1968})}\BibitemShut {NoStop}%
\bibitem [{\citenamefont {Gopalakrishnan}, \citenamefont {Kanevsky},\ and\
  \citenamefont {Nahamoo}(1989)}]{gopalakrishnan1989generalization}%
  \BibitemOpen
  \bibfield  {author} {\bibinfo {author} {\bibfnamefont {P.}~\bibnamefont
  {Gopalakrishnan}}, \bibinfo {author} {\bibfnamefont {D.}~\bibnamefont
  {Kanevsky}}, \ and\ \bibinfo {author} {\bibfnamefont {D.}~\bibnamefont
  {Nahamoo}},\ }in\ \href@noop {} {\emph {\bibinfo {booktitle} {Acoustics,
  Speech, and Signal Processing, 1989. ICASSP-89., 1989 International
  Conference on}}}\ (\bibinfo {year} {1989})\ p.\ \bibinfo {pages}
  {631}\BibitemShut {NoStop}%
\bibitem [{\citenamefont {Gangopadhyay}, \citenamefont {Chatterjee},\ and\
  \citenamefont {Chakrabartty}(2017)}]{svmgrowth}%
  \BibitemOpen
  \bibfield  {author} {\bibinfo {author} {\bibfnamefont {A.}~\bibnamefont
  {Gangopadhyay}}, \bibinfo {author} {\bibfnamefont {O.}~\bibnamefont
  {Chatterjee}}, \ and\ \bibinfo {author} {\bibfnamefont {S.}~\bibnamefont
  {Chakrabartty}},\ }\href@noop {} {\bibfield  {journal} {\bibinfo  {journal}
  {{I}EEE {T}ransactions on {N}eural {N}etworks and {L}earning {S}ystems}\ }
  (\bibinfo {year} {2017})},\ \bibinfo {note} {dOI:
  10.1109/TNNLS.2017.2690434}\BibitemShut {NoStop}%
\bibitem [{\citenamefont {Rockafellar}(1974)}]{rockafellar1974conjugate}%
  \BibitemOpen
  \bibfield  {author} {\bibinfo {author} {\bibfnamefont {R.~T.}\ \bibnamefont
  {Rockafellar}},\ }\href@noop {} {}\ (\bibinfo  {publisher} {SIAM},\ \bibinfo
  {year} {1974})\BibitemShut {NoStop}%
\end{thebibliography}%
\end{document}